\documentclass{article}
\usepackage{amsmath,amsfonts,amsthm,amssymb,amscd}

\renewcommand{\Re}{\mathop{\rm Re}\nolimits}

\newcommand{\p}{\partial}
\newcommand{\e}{\varepsilon}

\newcommand{\R}{{\mathbb R}}
\newcommand{\IP}{{\mathbb P}}
\newcommand{\Z}{{\mathbb Z}}
\newcommand{\E}{{\mathbb E}}

\newcommand{\I}{{\mathbb I}}

\newcommand{\FF}{{\cal F}}

\newcommand{\HH}{{\cal H}}

\theoremstyle{plain}
\newtheorem*{mt}{Main Theorem}
\newtheorem{theorem}{Theorem}[section]

\newtheorem{proposition}[theorem]{Proposition}
\newtheorem{corollary}[theorem]{Corollary}
\theoremstyle{definition}

\theoremstyle{remark}

\numberwithin{equation}{section}

\begin{document}
\author{Armen Shirikyan}
\title{Local times for solutions of the complex Ginzburg--Landau equation and the inviscid limit}
\date{\small D\'epartement de Math\'ematiques, Universit\'e de Cergy--Pontoise, CNRS UMR 8088\\
2 avenue Adolphe Chauvin, 95302 Cergy--Pontoise Cedex, France\\ 
E-mail: Armen.Shirikyan@u-cergy.fr}
\maketitle
\begin{abstract}
We consider the behaviour of the distribution for stationary solutions of the complex Ginzburg--Landau equation perturbed by a random force. It was proved in~\cite{KS-JPA2004} that if the random force is proportional to the square root of the viscosity~$\nu>0$, then the family of stationary measures possesses an accumulation point as~$\nu\to0^+$. We show that if~$\mu$ is such point, then the distributions of the~$L^2$ norm and of the energy possess a density with respect to the Lebesgue measure. The proofs are based on It\^o's formula and some properties of local time for semimartingales.


\end{abstract}



\section{Introduction}
\label{s1}
Let $D\subset\R^d$ be a bounded domain with a smooth boundary~$\p D$. In what follows, we always assume that $d\le 4$. 
Consider the complex Ginzburg--Landau equation in~$D$ perturbed by a white noise force smooth in the space variables: 
\begin{align}
\dot u-(\nu+i)\Delta u+i\lambda|u|^2u&=\sqrt\nu\,\eta(t,x),\label{1.1}\\
u\bigr|_{\p D}&=0.
\label{1.2}
\end{align}
Here $\nu\in(0,1]$ and $\lambda>0$ are some parameters, $u=u(t,x)$ is a complex-valued unknown function, and~$\eta$ is a random process of the form
\begin{equation} \label{1.3}
\eta(t,x)=\frac{\p}{\p t}\zeta(t,x), \quad
\zeta(t,x)=\sum_{j=1}^\infty b_j\beta_j(t)e_j(x),
\end{equation}
where $\{\beta_j=\beta_{j}^+ +i\beta_{j}^-\}$ is a sequence of independent complex-valued Brownian motions, 
$\{e_j\}$ is a complete set of normalised eigenfunctions of the Dirichlet Laplacian with eigenvalues~$\alpha_1<\alpha_2\le\alpha_3\le\cdots$, and~$b_j\ge0$ are some constants going to zero sufficiently fast.
Let $\{\mu_\nu\}$ be a family of stationary measures\footnote{The existence of a stationary measure for any $\nu>0$ follows from a priori estimates for solutions of~\eqref{1.1}, \eqref{1.2} and the Bogolyubov--Krylov argument; see~\cite{CK-1997,KS-JPA2004}.} for the Markov process associated with problem~\eqref{1.1}, \eqref{1.2}.  It was proved in~\cite{KS-JPA2004} that for any sequence $\hat\nu_k\to0^+$  the family~$\{\mu_{\hat\nu_k}\}$ has at least one accumulation point in the sense of weak convergence of measures on~$L^2(D)$, and any limiting point~$\mu$ for the family~$\{\mu_\nu\}$ satisfies the relations
\begin{align}
\mu(H^2)&=1,\label{1.5}\\
\int_{L^2}\|\nabla v\|^2\mu(dv)&=C_1,\label{1.6}\\
\int_{L^2}\bigl(\|v\|_{H^2}^2+\|v\|_{L^4}^4\bigr)\mu(dv)&\le C_2,
\label{1.7}
\end{align}
where $C_1$ and~$C_2$ are some constants depending on~$D$ and~$\{b_j\}$. Moreover, the measure~$\mu$ is the law of a stationary process $v(t,x)$ whose almost every realisation belongs to the space $L_{\rm loc}(\R_+,H^2)\cap W_{\rm loc}^{1,\frac43}(\R_+,L^{\frac43})$ and satisfies the nonlinear Schr\"odinger equation (NLS) obtained from~\eqref{1.1} by setting~$\nu=0$. Let us note that, in view of relation~\eqref{1.6} and the fact that the only steady state for NLS is the trivial zero solution, the support of~$\mu$ contains infinitely many points. In the case of the 2D Navier--Stokes equations, it was proved by Kuksin~\cite{kuksin-CMP2008} that if the random perturbation is non-degenerate, then the energy and enstrophy of such an inviscid limit have a density with respect to the Lebesgue measure. The present paper is devoted to the proof of similar results for~\eqref{1.1}, \eqref{1.2}. In particular, we establish the following theorem. 

\begin{mt}
Under the above hypotheses, the following assertions hold for any  family of stationary measures~$\{\mu_\nu\}$ and its limiting points~$\mu$ in the sense of weak convergence. 
\begin{itemize}
\item
If $b_j\ne0$ for some $j\ge1$, then the projection of~$\mu_\nu$ to the one-dimensional subspace spanned by the $j^{\text{th}}$ eigenfunction of the Direchlet Laplacian has a bounded density with respect to the Lebesgue measure. Moreover, the limiting measures~$\mu$ have no atom at $u=0$. 
\item
If $b_j\ne0$ for all $j\ge1$, then the laws of the functionals
\begin{align*}
\HH_0(u)&=\frac12\|u\|^2=\frac12\int_D|u(x)|^2dx,\\
\HH_1(u)&=\int_D\Bigl(\frac12\,|\nabla u(x)|^2
+\frac\lambda4\,|u(x)|^4\Bigr)\,dx
\end{align*}
under the probability distribution~$\mu$ on~$L^2$ possess a density with respect to the Lebesgue measure on the real line. 
\end{itemize}
\end{mt}

We refer the reader to Section~\ref{s3.1} for a more detailed statement of the results and to Sections~\ref{s3.2} and~\ref{s3.3} for the proofs. Note that the second assertion of the theorem was announced in Kuksin's paper~\cite{KS-JPA2004}. His proof for the case of the Navier--Stokes equations was based on a construction of an auxiliary stationary process satisfying a $\nu$-independent equation and an application of Krylov's estimate for semimartingales. Here we take a different approach using some properties of local times for the functionals in question. This enables one to simplify the proof and to get somewhat sharper estimates. The drawback of this approach is that it does not allow consideration of vector-valued functionals of solutions. 

In conclusion, let us mention that an alternative approach for constructing invariant measures for the Schr\"odinger equation with defocusing cubic nonlinearity was suggested by Tzvetkov~\cite{tzvetkov-2008}. His argument inspired by~\cite{LRS-1988} is based on considering renormalized Gibbs measures for finite-dimensional approximations and passing to the limit as the dimension goes to infinity. He proves that the resulting distribution is a non-degenerate Gaussian measure concentrated on radial functions of low regularity. 

\subsubsection*{Acknowledgements}
The author thanks L.~Zambotti for his remarks on the existence of density for martingales and useful references.  

\subsection*{Notation}
In what follows, we denote by~$D$ a bounded domain with smooth boundary, by~$\R$ the real axis, and by~$\Z_*$ the set of nonzero integers. We always assume that Polish spaces are endowed with their Borel $\sigma$-algebra and write~$\ell$ for the Lebesgue measure on~$\R$. Given any set $A$, we denote by~$\I_A$ its indicator function. We shall use the following functional spaces.

\smallskip
$L^q=L^q(D)$ stands for the space of complex-valued measurable functions~$u(x)$ such that 
$$
\|u\|_{L^q}=\biggl(\int_D|u(x)|^qdx\biggr)^{1/q}<\infty.
$$
We regard~$L^2$ as a real Hilbert space with the scalar product 
$$
(u,v)=\Re\int_D u(x)\overline{v(x)}\,dx
$$
and denote by~$\|\cdot\|$ the corresponding  norm.

\smallskip
$H^s=H^s(D)$ denotes the Sobolev space of order~$s$ endowed with the usual norm~$\|\cdot\|_s$. 

\smallskip
$\{e_j,j\ge1\}$ stands for the complete set of normalised eigenfunctions of the Dirichlet Laplacian. We denote $e_{-j}=ie_j$ for $j\le-1$, so that $\{e_j,j\in\Z_*\}$ is an orthonormal basis in~$L^2$. 

\section{Preliminaries}
\label{s2}

\subsection{Stationary measures and a priori estimates}
\label{s2.1}
In this subsection, we have collected some known facts about the existence of stationary solutions for problem~\eqref{1.1}, \eqref{1.2} and a priori estimates for them. We shall assume that the coefficients~$b_j\ge0$ entering~\eqref{1.3} satisfy the inequalities\,\footnote{Of course, the first inequality is a consequence of the second. We wrote it to define the constant~$B_0$.}
\begin{equation} \label{2.1}
B_0=\sum_{j=1}^\infty b_j^2<\infty, \quad 
B_1=\sum_{j=1}^\infty \alpha_j b_j^2<\infty,\quad 
M=\sup_{x\in D}\sum_{j=1}^\infty b_j^2e_j^2(x)<\infty.
\end{equation}
In this case, it was proved in~\cite{KS-JPA2004} that the Cauchy problem~\eqref{1.1}, \eqref{1.2} is well posed in the space~$H_0^1$, and the corresponding Markov family has at least one stationary distribution for any $\nu>0$. The results listed in the following proposition are either established in the papers~\cite{KS-JPA2004,odasso-2006,shirikyan-dcds2006} or can be proved with the help of the methods used there. 

\begin{proposition} \label{p4.0}
Under the above hypotheses, for any stationary solution~$u(t)$ of~\eqref{1.1}, \eqref{1.2} the random processes~$\HH_0(u)$ and~$\HH_1(u)$ are semimartingales and can be represented in the form
\begin{align}
&\HH_0(u(t))=\HH_0(u(0))+\nu\int_0^t\bigl(B_0-\|\nabla u(t)\|^2\bigr)ds
+\sqrt\nu\sum_{j\in\Z_*}b_j\int_0^tu_jdw_j,\label{1.06}\\
&\HH_1(u(t))=\HH_1(u(0))+\nu\int_0^t\biggl(B_1-\|\Delta u\|^2
-2\lambda\bigl(|u|^2,|\nabla u|^2\bigr)
-\lambda\bigl(u^2,(\nabla u)^2\bigr)
\notag\\
&\quad+2\lambda\sum_{j=1}^\infty
b_j^2\bigl(|u|^2,e_j^2\bigr)\biggr)\,ds
+\sqrt\nu\sum_{j\in\Z_*}\int_0^tb_j\bigl(-\Delta u+\lambda|u|^2u,e_j\bigr)dw_j,
\label{1.07}
\end{align}
where we set $b_{-j}=b_j$ for $j\ge1$, $u_j=(u,e_j)$, and $w_j=\beta_{j}^\pm$ for $\pm j\ge1$. 
Moreover, there is a constant~$C>0$ not depending~$\nu$ such that 
\begin{align}
\E\,\|\nabla u\|^2&=B_0,\label{2.41}\\
\E\,\bigl(\|u\|_{L^6}^3+\bigl(|u|^2,|\nabla u|^2\bigr)+\|u\|_2^2\bigr)
&\le C (B_1+MB_0), \label{2.42}
\end{align}
\end{proposition}

\subsection{Local time for semimartingales}
\label{s2.2}
Let $(\Omega,\FF,\IP)$ be a complete probability space with a right-continuous filtration~$\{\FF_t, t\ge0\}$ augmented with respect to~$(\FF,\IP)$, let~$\{\beta_j\}$ be a sequence of independent Brownian motions with respect to~$\FF_t$, and let~$y_t$ be a scalar semimartingale of the form
\begin{equation} \label{7.36}
y_t=y_0+\int_0^tx_sds+\sum_{j=1}^\infty\int_0^t\theta_s^jd\beta_j,
\end{equation}
where $x_t$ and $\theta_t^j$ are $\FF_t$-adapted processes such that 
$$
\E\int_0^t\biggl(x_s^2+\sum_{j=1}^\infty\bigl|\theta_s^j\bigr|^2\biggr)\,ds<\infty
\quad\mbox{for any $t>0$}.
$$
The following result is a straightforward consequence of Theorem~7.1 in~\cite[Chapter~3]{KS1991}. 

\begin{theorem} \label{t7.6.1}
Under the above hypotheses, then the random field 
\begin{equation} \label{7.49}
 \Lambda_t(a,\omega)=|y_t-a|-|y_0-a|-\sum_{j=1}^\infty\int_0^t\I_{(a,\infty)}(y_s)\theta_s^jd\beta_j-\int_0^t\I_{(a,\infty)}(y_s)x_sds
\end{equation}
defined for $t\ge0$ and $a\in\R$ possesses the following properties hold. 
\begin{itemize}
\item[\bf(i)]
The mapping $(t,a,\omega)\mapsto\lambda_t(a,\omega)$ is measurable, and for any $a\in\R$ the process $t\mapsto\Lambda_t(a,\omega)$ is $\FF_t$-adapted, continuous and non-decreasing. Moreover, for every $t\ge0$ and almost every~$\omega\in\Omega$ the function $a\mapsto\Lambda_t(a,\omega)$ is right-continuous.  
\item[\bf(ii)]
For any non-negative Borel-measurable function $h:\R\to\R$, with probability~$1$ we have 
\begin{equation} \label{7.48}
2\int_{-\infty}^\infty h(a)\Lambda_t(a,\omega)da=\sum_{j=1}^\infty\int_0^t h(y_s)|\theta_s^j|^2ds, \quad t\ge0.
\end{equation}
\end{itemize}
\end{theorem}
The random field $\Lambda_t(a,\omega)$ is called a  {\it local time\/} for~$y_t$, and~\eqref{7.49} is usually referred to as the {\it change of variable formula\/}.

\section{Main results}
\label{s3}

\subsection{Formulation}
\label{s3.1}

The following theorem, which is the main result of this paper, describes some qualitative properties of stationary measures for~\eqref{1.1}. 

\begin{theorem} \label{t3.1}
Assume that the coefficients~$b_j$ entering the definition of the random process~$\zeta$ satisfy inequalities~\eqref{2.1}. Then the following assertions hold for any stationary measure~$\mu_\nu$ of problem~\eqref{1.1}, \eqref{1.2} with $\nu>0$. 
\begin{itemize}
\item[\bf(i)]
Let $b_k\ne0$ for some $k\ge1$ and let $v\in L^2$ be a function non-orthogonal to~$e_j$. Then
the projection of~$\mu_\nu$ to the vector space spanned by~$v$ has a bounded density with respect to the Lebesgue measure. In particular, $\mu_\nu$~has no atoms for $\nu>0$. 
\item[\bf(ii)]
Let $b_k\ne0$ for some $k\ge1$. Then there is constant $C>0$ not depending on the sequence~$\{b_j\}$ and the parameter $\nu>0$ such that 
\begin{equation} \label{3.0}
\mu_\nu\bigl(\{u\in L^2:\|u\|_{L^2}\le \delta\}\bigr)\le CB_0^{-1} \sqrt{B_1+MB_0}\,\delta\quad\mbox{for any}\quad\delta\ge0.
\end{equation}
\item[\bf(iii)]
Let $b_j\ne0$ for all $j\ge1$. Then there is a continuous increasing function~$p(r)$ going to zero with~$r$ such that
\begin{equation} \label{2.2}
\mu_\nu\bigl(\{u\in L^2:\HH_0(u)\in\Gamma\}\bigr)+
\mu_\nu\bigl(\{u\in L^2:\HH_1(u)\in\Gamma\}\bigr)\le p\bigl(\ell(\Gamma)\bigr)
\end{equation}
for any Borel subset $\Gamma\subset\R$.
\end{itemize}
\end{theorem}

Let us emphasize that the stationarity of~$u$ is important for the validity of the conclusions of Theorem~\ref{t3.1}. A counterexample constructed by Fabes and Kenig~\cite{FK-1981} shows that a solution of a one-dimensional SDE with a diffusion term separated from zero may have a singular distribution (see also~\cite{martini-2000}).
 
 \smallskip
Theorem~\ref{t3.1} immediately implies the results formulated in the Introduction. Moreover, we have the following corollary. 

\begin{corollary} \label{c3.3}
Under the hypotheses of Theorem~\ref{t3.1}, the following assertions hold for any limiting point~$\mu$ of the family~$\{\mu_\nu,\nu>0\}$  in the sense of weak convergence on~$L^2$. 
\begin{itemize}
\item
Let $b_k\ne0$ for some $k\ge1$. Then~$\mu$ has no atom at $u=0$ and satisfies inequality~\eqref{3.0} in which~$\mu_\nu$ is replaced by~$\mu$.  
\item 
Let $b_j\ne0$ for all $j\ge1$. Then for any Borel subset $\Gamma\subset\R$ we have
\begin{equation} \label{3.3}
\mu\bigl(\{u\in L^2:\HH_0(u)\in\Gamma\}\bigr)+
\mu\bigl(\{u\in L^2:\HH_1(u)\in\Gamma\}\bigr)\le p\bigl(\ell(\Gamma)\bigr),
\end{equation}
where~$p$ is the function constructed in Theorem~\ref{t3.1}. 
\end{itemize} 
\end{corollary}

\begin{proof}
We shall confine ourselves to the proof of the second assertion, because the first one can be established by a similar argument. 
It is well known that if a sequence of measures~$\mu_{\nu_k}$ converges to~$\mu$ weakly on~$L^2$, then 
$$
\liminf_{k\to\infty}\mu_{\nu_k}(G)\le\mu(G)\quad
\mbox{for any open subset $G\subset L^2$}. 
$$
Combining this with~\eqref{3.0}, we see that~\eqref{3.3} holds for any open subset $\Gamma\subset\R$. Now recall that if~$\lambda$ is a Borel measure on~$\R$, then 
$$
\lambda(\Gamma)=\inf\{\lambda(G): G\supset\Gamma, \mbox{$G$ is open}\}.
$$
Combining this property with the continuity of~$p$, we conclude that~\eqref{3.3} is true for any Borel set~$\Gamma\subset\R$. 
\end{proof}

The proof of Theorem~\ref{t3.1} is given in Section~\ref{s3.3}. A key ingredient of the proof is the following result established in the next subsection.

\begin{proposition} \label{p3.3}
Let us assume that the hypotheses of Theorem~\ref{t3.1} are satisfied and let $g\in C^2(\R)$ be a real-valued function of at most polynomial growth at infinity. Then for any stationary solution~$u(t,x)$ for~\eqref{1.1}, \eqref{1.2} and any Borel subset $\Gamma\subset\R$ we have
\begin{multline} \label{2.3}
\E\int_\Gamma\I_{(a,\infty)}\bigl(g(\|u\|^2)\bigr)
\biggl(g'(\|u\|^2)(B_0-\|\nabla u\|^2)+g''(\|u\|^2)\sum_{j\in\Z_*}b_j^2u_j^2\biggr)\,da\\
+\sum_{j\in\Z_*}b_j^2\,\E\Bigl(\I_\Gamma\bigl(g(\|u\|^2)\bigr)\bigl(g'(\|u\|^2)\,u_j\bigr)^2\Bigr)=0.
\end{multline}
\end{proposition}

\subsection{Proof of Proposition~\ref{p3.3}}
\label{s3.2}
Let us fix any function $g\in C^2(\R)$ and consider the process $f(t)=g(u(t))$, where $u(t)$ is a stationary solution of~\eqref{1.1}, \eqref{1.2}. It follows from~\eqref{1.06} that~$f$ is a semimartingale which can be written as
\begin{equation} \label{3.4}
f(t)=f(0)+\nu\int_0^tA(s)\,ds+2\sqrt\nu\sum_{j\in\Z_*}b_j\int_0^tg'(\|u\|^2)u_jdw_j,
\end{equation}
where we set
$$
A(t)=2\biggl(g'(\|u\|^2)\bigl(B_0-\|\nabla u\|^2\bigr)+g''(\|u\|^2)\sum_{j\in\Z_*}b_j^2u_j^2\biggr).
$$
Let $\Lambda_t(a)$ be the local time for~$f$. Then, in view of relation~\eqref{7.48} with~$h=\I_\Gamma$, we have
\begin{equation} \label{3.5}
2\int_\Gamma \Lambda_t(a)\,da=4\nu\sum_{j\in\Z_*}b_j^2\int_0^t\I_\Gamma(f(s))\bigl(g'(\|u\|^2)u_j\bigr)^2\,ds. 
\end{equation}
Taking the mean value and using the stationarity of~$u$, we derive
\begin{equation} \label{3.6}
\int_\Gamma \bigl(\E\Lambda_t(a)\bigr)\,da
=2\nu t\sum_{j\in\Z_*} b_j^2 \,\E\Bigl(\I_\Gamma(f)\bigl(g'(\|u\|^2)u_j\bigr)^2\Bigr). 
\end{equation}
On the other hand, by the change of variable formula~\eqref{7.49}, for any $a\in\R$ we have
\begin{align*}
\Lambda_t(a)
&=|f(t)-a|-|f(0)-a|-2\sqrt\nu \sum_{j\in\Z_*} b_j \int_0^t\I_{(a,\infty)}(f(s))g'(\|u\|^2)u_jdw_j\\
&\qquad -\nu\int_0^t \I_{(a,\infty)}(f(s))A(s)\,ds.
\end{align*}
Taking the mean value and using again the stationarity of~$u$, we derive
$$
\E \Lambda_t(a)
=-\nu t\,\E\bigl(\I_{(a,\infty)}(f(0))\,A(0))\bigr).
$$
Substituting this into~\eqref{3.6} and recalling the definition of~$A$, we arrive at the required relation~\eqref{2.3}.

\subsection{Proof of Theorem~\ref{t3.1}}
\label{s3.3}
{\it Proof of~{\rm(i)}}. 
We repeat essentially the argument used in the proof of Proposition~\ref{p3.3}. Let $v\in  L^2$ be a function non-orthogonal to~$e_k$ or $e_{-k}$ and let~$\mu_\nu^v$ be the projection of~$\mu_\nu$ to the vector space spanned by~$v$. Then, in view of~\eqref{1.1}, the stationary process $z(t)=(u(t),v)$ is a semimartingale, its law coincides with~$\mu_\nu^v$, and it can be written in the form
$$
z(t)=z(0)+\int_0^tg(s)\,ds+\sqrt\nu\sum_{j\in\Z_*}d_jdw_j,
$$
where we set
$$
g(t)=\bigl((\nu+i)\Delta u-\lambda |u|^2u,v\bigr), \quad d_j=b_j(v,e_j).
$$
Note that either $d_k\ne0$ or $d_{-k}\ne0$, whence it follows that 
$$
D_0:=\sum_{j\in\Z_*} d_j^2\ne0. 
$$
Let $\Lambda_t^v(a)$ be the local time for~$z$. Then, by~\eqref{7.48}, for any Borel subset $\Gamma\subset\R$, we have
$$
2\int_{\Gamma} \Lambda_t^v(a)\,da
=\sqrt\nu\, D_0\int_0^t I_\Gamma(z(s))\,ds, \quad t\ge0.
$$
Taking the mean value and using the stationarity of~$z$, we see that
\begin{equation} \label{3.7}
2\int_{\Gamma} \bigl(\E\Lambda_t^v(a)\bigr)\,da
=\sqrt\nu\, t D_0 \mu_\nu^v(\Gamma), \quad t\ge0.
\end{equation}
On the other hand, by the change of variable formula~\eqref{7.49}, we have
\begin{align*} 
 \Lambda_t^v(a)&=|z(t)-a|-|z(0)-a|-\int_0^t\I_{(a,\infty)}(z(s))g(s)\,ds\\
&\qquad -\sqrt\nu\sum_{j\in\Z_*} d_j\int_0^t\I_{(a,\infty)}(z(s))dw_j. 
\end{align*}
Taking the mean value and substituting the resulting formula into~\eqref{3.7}, we derive
\begin{equation} \label{3.9}
\IP\bigl(\{z(0)\in\Gamma\}\bigr)
\le\frac{2}{\sqrt\nu D_0}\bigl(\E\,|g(0)|\bigr)\ell(\Gamma). 
\end{equation}
Using~\eqref{2.42}, we see that
\begin{equation} \label{3.10}
\E\,|g(0)|\le C_1\E\bigl(\|u\|_{L^6}^3+\|u\|_2^2\bigr)\le C_2.
\end{equation}
Since the law of~$z(0)$ coincides with~$\mu_\nu^v$, combining~\eqref{3.9} and~\eqref{3.10}, we conclude that
$$
 \mu_\nu^v(\Gamma)\le C_3\nu^{-1/2}\ell(\Gamma)\quad\mbox{for any Borel set $\Gamma\subset\R$}. 
$$
This inequality implies that $\mu_\nu^v$ has a bounded density with respect to the Lebesgue measure. 

\medskip
{\it Proof of~{\rm(ii)}}. 
Let us apply relation~\eqref{2.3} in which $\Gamma=[\alpha,\beta]$ with $\alpha>0$ and $g\in C^2(\R)$ is a function that coincides with~$\sqrt{x}$ for~$x\ge\alpha$. This results in
\begin{multline*}
\E\int_\alpha^\beta \I_{(a,\infty)}\bigl(\|u\|\bigr)
\biggl(\frac{B_0-\|\nabla u\|^2}{2\|u\|}-\frac{1}{4\|u\|^3}\sum_{j\in\Z_*}b_j^2u_j^2\biggr)\,da\\
+\frac14\sum_{j\in\Z_*}b_j^2\,\E\Bigl(\I_\Gamma\bigl(\|u\|\bigr)\|u\|^{-2}u_j^2 \Bigr)=0.
\end{multline*}
It follows that
\begin{equation} \label{3.11}
\E\int_\alpha^\beta\frac{\I_{(a,\infty)}(\|u\|)}{\|u\|^3}
\biggl(2B_0\|u\|^2-\sum_{j\in\Z_*}b_j^2u_j^2\biggr)\,da
\le 2(\beta-\alpha)\,\E\biggl(\frac{\|\nabla u\|^2}{\|u\|}\biggr).
\end{equation}
Now note that
\begin{align*}
2B_0\|u\|^2-\sum_{j\in\Z_*}b_j^2u_j^2&=\sum_{j\in\Z_*}(2B_0-b_j^2)u_j^2\ge B_0\|u\|^2,\\
\E\biggl(\frac{\|\nabla u\|^2}{\|u\|}\biggr)&\le C_4\E\,\|\Delta u\| \le C_5 \sqrt{B_1+MB_0}, 
\end{align*}
where we used interpolation and inequality~\eqref{2.42}. Substituting these estimates into~\eqref{3.11} and passing to the limit as $\alpha\to0^+$, we derive
\begin{equation} \label{3.12}
\E \int_0^\beta \I_{(a,\infty)}(\|u\|) \|u\|^{-1}da\le C_6 B_0^{-1} \sqrt{B_1+MB_0}\,\beta. 
\end{equation}
We now fix a constant~$\delta>0$ and note that the left-hand side of~\eqref{3.12} can be minorised by
\begin{align*}
\E \int_0^\beta \I_{(a,\delta]}(\|u\|) \|u\|^{-1}da
&\ge \delta^{-1}\E \int_0^\beta \I_{(a,\delta]}(\|u\|)\,da \\
&= \delta^{-1} \int_0^\beta \IP\bigl(\{a<\|u\|\le\delta)\,da. 
\end{align*}
Substituting this inequality into~\eqref{3.12}, we obtain
$$
\frac{1}{\beta} \int_0^\beta \IP\bigl(\{a<\|u\|\le\delta)\,da
\le C_6 B_0^{-1} \sqrt{B_1+MB_0}\,\delta. 
$$
Passing to the limit as $\beta\to 0^+$ and recalling that~$\mu_\nu$ has no atom at $u=0$, we arrive at the required inequality~\eqref{3.0}. 

\medskip
{\it Proof of~{\rm(iii)}}. 
It suffices to show that each term on the left-hand side of~\eqref{2.2} can be estimated by~$p(\ell(\Gamma))$. We begin with the case of the functional $\HH_0(u)$. 

Applying relation~\eqref{2.3} with $g(x)=x$ and~$\Gamma$ replaced by$~2\Gamma$, we obtain
\begin{equation} \label{3.13}
\E\biggl(\I_\Gamma\bigl(\HH_0(u)\bigr)\sum_{j\in\Z_*}b_j^2u_j^2\biggr)
\le \int_{2\Gamma}\E\Bigl(\I_{(a,\infty)}\bigl(\|u\|^2\bigr)\|\nabla u\|^2\Bigr)\,da\le 2B_0\ell(\Gamma), 
\end{equation}
where we used~\eqref{2.41} to get the second inequality. We wish to estimate the left hand side of this inequality from below. To this end, note that if $\|u\|\ge\delta$ and $\|u\|_2\le\delta^{-1/2}$, then for any integer $N\ge1$ we have
\begin{align*}
\sum_{j\in\Z_*}b_j^2u_j^2
&\ge \underline{b}_N^2\sum_{0<|j|\le N}u_j^2
=\underline{b}_N^2\biggl(\|u\|^2-\sum_{|j|>N}u_j^2\biggr)\\
&\ge \underline{b}_N^2\biggl(\|u\|^2-\alpha_{N+1}^{-2}\|\Delta u\|^2\biggr)
\ge \underline{b}_N^2\bigl(\delta^2-\alpha_{N+1}^{-2}\delta^{-1}\bigr),
\end{align*}
where $\underline{b}_N=\min\{b_j,1\le j\le N\}$. Choosing $N=N(\delta)$ sufficiently large, we find an increasing function~$\e(\delta)>0$ going to zero with~$\delta$ such that 
\begin{equation} \label{3.14}
\sum_{j\in\Z_*}b_j^2u_j^2\ge \e(\delta)\quad
\mbox{for $\|u\|\ge\delta$, $\|u\|_2\le\delta^{-1/2}$}.
\end{equation}

Define now the event $G_\delta=\{\|u(0)\|\le\delta\mbox{ or } \|u(0)\|_2\ge\sqrt\delta\}$ and note that, in view of~\eqref{3.0}, \eqref{2.42}, and Chebyshev's inequality, we have
$$
\IP(G_\delta)\le \IP\{\|u\|\le\delta\}+\IP\{\|u\|\ge\sqrt\delta\}\le C_7\delta.
$$
Combining this with~\eqref{3.13} and~\eqref{3.14}, we write
\begin{align*}
\IP\{\HH_0(u)\in\Gamma\}
&=\IP\bigl(\{\HH_0(u)\in\Gamma\}\cap G_\delta\bigr)
+\IP\bigl(\{\HH_0(u)\in\Gamma\}\cap G_\delta^c\bigr)\\
&\le C_7\delta +\e(\delta)^{-1}\E\biggl(\I_\Gamma\bigl(\HH_0(u)\bigr)\sum_{j\in\Z_*}b_j^2u_j^2\biggr)\\
&\le C_7\delta +C_8\e(\delta)^{-1}\ell(\Gamma). 
\end{align*}
This inequality immediately implies the required result. 

\medskip
To prove the estimate for~$\HH_1(u)$, we use a similar argument; however, the calculations become more involved. A literal repetition of the proof of Proposition~\ref{p3.3} enables one to show that 
\begin{align} 
\sum_{j\in\Z_*}b_j^2\,\E\Bigl(\I_\Gamma\bigl(\HH_1(u)\bigr)\bigl(-\Delta u+\lambda|u|^2u,e_j\bigr)^2\Bigr)
&\le C_8\ell(\Gamma)\E\Bigl(\|\Delta u\|^2+\bigl(|u|^2,|\nabla u|^2\bigr)\Bigr)\notag\\
&\le C_9 \bigl(B_1+MB_0\bigr)\ell(\Gamma), \label{3.15} 
\end{align}
where we used~\eqref{2.42} to derive the second inequality. Let us estimate from below the following expression arising in the left-hand side of~\eqref{3.15}:
$$
\Xi(u)=\sum_{j\in\Z_*}b_j^2\bigl(-\Delta u+\lambda|u|^2u,e_j\bigr)^2.
$$
Denoting by~$A_u$ the operator $-\Delta+\lambda|u|^2$, we see that $(A_uv,v)\ge\alpha_1\|v\|^2$ for any $v\in H^2\cap H_0^1$. It follows that
\begin{equation} \label{3.16}
\sum_{j\in\Z_*}(A_uu,e_j)u_j\ge \alpha_1\|u\|^2. 
\end{equation}
On the other hand, for any integer~$N\ge1$, we can write
\begin{align} 
\sum_{j\in\Z_*}(A_uu,e_j)u_j
&=\sum_{0<|j|\le N}(A_uu,e_j)u_j+(A_uu,Q_Nu)\notag\\
&\le{\underline b}_N^{-1}\|u\| \biggl(\,\sum_{0<|j|\le N}b_j^2(A_uu,e_j)^2\biggr)^{1/2}+\|A_uu\|\,\|Q_Nu\|\notag\\
&\le{\underline b}_N^{-1}\|u\|\sqrt{\Xi(u)}+\alpha_N^{-1}\|\Delta u\|\bigl(\|\Delta u\| +\lambda\|u\|_{L^6}^3\bigr),
\label{3.17}
\end{align}
where $Q_N:L^2\to L^2$ stands for the orthogonal projection onto the subspace spanned by~$e_j$, $|j|>N$, and we used the inequality $\|Q_Nu\|\le\alpha_N^{-1}\|\Delta u\|$. Combining~\eqref{3.16} and~\eqref{3.17}, we see that, if $\|u\|\ge\delta$ and $\|\Delta u\|+\|u\|_{L^6}^3\le\delta^{-1}$, then 
\begin{align*}
\sqrt{\Xi(u)}
&\ge {\underline b}_N\Bigl(\alpha_1\|u\|-\alpha_N^{-1}\|u\|^{-1}\|\Delta u\|\bigl(\|\Delta u\|+\lambda\|u\|_{L^6}^3\bigr)\Bigr)\\
&\ge {\underline b}_N\bigl(\alpha_1\delta-\alpha_N^{-1}(1+\lambda)\delta^{-3}\bigr).
\end{align*}
Choosing $N=N(\delta)$ sufficiently large, we find an increasing function~$\e(\delta)>0$ going to zero with~$\delta>0$ such that
$$
\sum_{j\in\Z_*}b_j^2\bigl(-\Delta u+\lambda|u|^2u,e_j\bigr)^2\ge\e(\delta)\quad
\mbox{for $\|u\|\ge\delta$, $\|\Delta u\|+\|u\|_{L^6}^3\le\delta^{-1}$}. 
$$
The required upper bound for $\IP\{\HH_1(u)\in\Gamma\}$ can now be derived from~\eqref{3.15} in exactly the same way as for the case of~$\HH_0$. The proof of Theorem~\ref{t3.1} is complete. 

\def\cprime{$'$} \def\cprime{$'$}
\providecommand{\bysame}{\leavevmode\hbox to3em{\hrulefill}\thinspace}
\providecommand{\MR}{\relax\ifhmode\unskip\space\fi MR }

\providecommand{\MRhref}[2]{%
  \href{http://www.ams.org/mathscinet-getitem?mr=#1}{#2}
}
\providecommand{\href}[2]{#2}

\end{document}